\documentclass[12pt]{article}

\usepackage{amsmath,amsfonts}
\usepackage{graphicx}

\newcommand{\CC}{{\cal{C}}}

\newcommand{\R}{\mathbb{R}}

\renewcommand{\SS}{{\cal{S}}}

\def\version{1.0}
\def\versiondate{July 1, 2008}
 
\title{POCP: a package for\\
Polynomial Optimal Control Problems}
\date{User's guide. Version \version~of \versiondate }
\author{Didier Henrion, Jean B. Lasserre and Carlo Savorgnan
\thanks{D. Henrion is with LAAS-CNRS, University of Toulouse, France
and also with the Faculty of Electrical Engineering,
Czech Technical University in Prague, Czech Republic.
J. B. Lasserre is with LAAS-CNRS, University of Toulouse, France and 
also with the Institute of Mathematics of the University of Toulouse, France.
C. Savorgnan was with LAAS-CNRS, University of Toulouse, France.
He is now with the Department of Electrical Engineering,
Katholieke Universiteit Leuven, Belgium.
This work was partly funded by project MOGA of
the French National Research Agency (ANR) and
by project No.~102/06/0652
of the Grant Agency of the Czech Republic.}}

\begin{document}
\maketitle

\begin{abstract}
POCP is a new Matlab package running jointly with GloptiPoly 3 and, optionally, YALMIP. 
It is aimed at nonlinear optimal control problems for which all the problem data are polynomial,
and provides an approximation of the optimal value as well as some control policy. 
Thanks to a user-friendly interface, POCP reformulates such control problems as generalized problems
of moments, in turn converted by GloptiPoly 3 into a hierarchy of 
semidefinite programming problems whose associated sequence of optimal values
converges to the optimal value of the polynomial optimal control problem.
In this paper we describe the basic features of POCP and illustrate them with some
numerical examples.
\end{abstract}

\section{What is POCP?}
POCP is a Matlab package aimed at solving approximately
nonlinear optimal control problems (OCPs)
for which all the problem data are polynomial.
Consider a continuous-time system described by the differential equation
\begin{equation}\label{eq:dynamics}
\dot x(t) = f(t, x(t), u(t)),
\end{equation}
where $x\in\R^n$ and $u\in\R^m$ are the state and input vectors, respectively.
Given the running cost $h:\R\times\R^n\times\R^m \rightarrow \R$ and the final cost $H:\R^{n}\rightarrow\R$, the total 
cost to be minimized in the OCP is defined as
\begin{equation}\label{eq:cost}
J(0,T,x(0),u(\cdot))=\int_0^T h(t, x(t), u(t)) dt + H(x(T)),
\end{equation}
where $T$ is the horizon length.
Different OCPs can be formulated depending if the initial condition $x(0)$ has been assigned or not.

\subsection{Optimal control problem with free initial condition}

Consider the following constraints:
\begin{itemize}
\item $x(0) \in \CC_I=\{ x: g_{I_j}(x)\leq 0, ~ j=1,\dots,n_I \} $;
\item $x(T) \in \CC_F=\{ x: g_{F_j}(x)\leq 0, ~ j=1,\dots,n_F \}$;
\item $(t, x(t), u(t)) \in \CC_T=\{ (t,x,u): g_{T_j}(t,x,u)\leq 0, ~ j=1,\dots,n_T \}$.
\end{itemize}
An OCP with free initial condition is of the form
\begin{equation}\label{eq:probfree}
\begin{array}{ll}
\displaystyle\min_{x(0), u(t)} & J(0,T,x(0),u(t)) \\
& x(0)\in \CC_{I} \\
& (t,x(t),u(t))\in \CC_{T} \\
& x(T)\in \CC_{F}
\end{array}
\end{equation}
To solve this problem using POCP all the data must be polynomial.
More specifically, $f$, $h$, $H$, $g_{I_j}$, $g_{T_j}$, and $g_{F_j}$ must be polynomial functions.

\subsection{Optimal control problem with fixed initial condition}

When set $\CC_{I}$ contains only one point, the initial condition is fixed. From this viewpoint, the minimization problem (\ref{eq:probfree}) can be used also to formulate OCPs where the initial condition is fixed.

However, it can be interesting to consider a more general framework where the initial condition
is not known exactly, but stochastically. In this case, the probability measure of the initial condition $\mu_{I}$ is given,
and we consider the following OCP 
\begin{equation}\label{eq:probfixed}
\begin{array}{ll}
\displaystyle\min_{u(\cdot)} & \int J(0,T,x_{0},u(x_{0},t)) d\mu_{I}(x_{0}) \\
& (t,x(t),u(t))\in \CC_{T} \\
& x(T)\in \CC_{F}.
\end{array},
\end{equation}
To better understand how to use this formulation, see \cite{HenLasSav2008}.

POCP can also deal with problems where some state variables are fixed at time $t=0$ and some of them can be chosen inside a set $\CC_I$.

POCP can be used for problems where the horizon $T$ is fixed or not.

\section{How does POCP work and what does it calculate?}
A full explanation of the theory on which POCP is based is out of scope for this user's guide.
The interested reader is referred to \cite{HenLasSav2008} and \cite{LHPT2008}.

The techniques implemented in POCP are based on a modeling via
occupation measures associated with semidefinite programming (SDP) relaxations. POCP formulates a hierarchy of relaxations giving a sequence of lower bounds on the optimal value of the OCPs (\ref{eq:probfree}) and (\ref{eq:probfixed}). Furthermore, POCP returns the moment matrices of the occupation measures used and, as a byproduct, a polynomial subsolution of the Hamilton-Jacobi-Bellman (HJB) equation.

\section{Installation}
POCP is a free Matlab package consisting of an archive file
downloadable from
\begin{center}
\begin{verbatim}
www.laas.fr/~henrion/software/pocp
\end{verbatim}
\end{center}
The installation consists of two steps:
\begin{itemize}
\item extract the directory \texttt{@pocp} from the archive file,
\item copy \texttt{@pocp} on your working directory or, using the command \texttt{addpath}, add the parent directory of \texttt{@pocp} to the Matlab path.
\end{itemize}
POCP is based on GloptiPoly version 3.3 or higher \cite{HenLasLof2007}.
Therefore it is assumed that this package is properly installed.
To compute polynomial subsolutions of the HJB equation,
the optimization modelling toolbox YALMIP \cite{Lof2004}
must be installed as well.

\section{An introductory example}
Consider the double integrator
\begin{equation}\nonumber
\begin{bmatrix} \dot x_1(t) \\ \dot x_2(t) \end{bmatrix} =
\begin{bmatrix} 0 & 1 \\ 0 & 0 \end{bmatrix}
\begin{bmatrix} x_1(t) \\ x_2(t) \end{bmatrix} +
\begin{bmatrix} 0 \\ 1 \end{bmatrix} u(t)
\end{equation}
with state constraint $x_2(t)\geq -1$ and input constraint $-1\leq u(t) \leq 1$.
When we consider the problem of steering in minimum time the state to the origin, an analytic solution to this problem is available \cite{LHPT2008}. When $x(0)=[1~1]^T$, the value of the minimum time to reach the origin is $T=3.5$.
In the following, we show a Matlab POCP script which defines the problem and calculates a lower bound on the optimal value.

\begin{verbatim}
% system variables
mpol x 2;  % state vector of size 2
mpol u; % scalar input

% state and input matrices
A = [0 1; 0 0];
B = [0; 1];

% define optimal control problem
prob = pocp( ...
       'state', x, ...
       'input', u, ...
       'dynamics', A*x+B*u);

% set further properties
prob = set(prob, 'idirac', x, [1; 1], ... % initial condition
                 'fdirac', x, [0; 0], ... % final condition
                 'tconstraint', [x(2)>=-1; u>=-1; u<=1], ... % constraints 
                 'scost', 1);  % setting integral cost h to 1

% call to the solver
[status, cost] = solvepocp(prob, 14); % 14 = degree of moment matrix

disp('The lower bound on the minimum time is');
disp(cost);
\end{verbatim}
We can notice that the user is required to know only four commands: \texttt{mpol} (a GloptiPoly command to define  variables), \texttt{pocp}, \texttt{set} and \texttt{solvepocp}.
When the script file is executed, we obtain the following output
\begin{verbatim}
...
The lower bound on the minimum time is
    3.4988
\end{verbatim}
In the following section we illustrate all POCP commands together with examples of their usage.

\section{Command guide}
POCP consists of three commands: \texttt{pocp}, \texttt{set}, and \texttt{solvepocp}.

\subsection{The \texttt{pocp} command}
\textbf{Short description:} Creates a new POCP object.
\\
\textbf{Syntax and usage:} When used in the form
\begin{verbatim}
>> prob = pocp;
\end{verbatim}
the command
\texttt{pocp} creates an empty POCP object \texttt{prob} of class \texttt{pocp}.
The problem data can be then specified using the command \texttt{set}.

The \texttt{pocp} command can also be used to create a POCP object and, at the same time, set the problem data. In this case, the syntax is
\begin{verbatim}
>> prob = pocp('Property1', Argument1, Argument2, ...,'Property2', ...)
\end{verbatim}
The properties and the arguments should be specified accordingly to the syntax of the \texttt{set} command,
see below.

\subsection{The \texttt{set} command}
\textbf{Short description:} Sets the data of a POCP.
\\
\textbf{Syntax and usage:} Specifies one or more properties of the POCP object. For example,
to set a property which requires two arguments, use
\begin{verbatim}
>> prob = set(prob, 'Property1', Argument1, Argument2);
\end{verbatim}
To set more properties at once,
\begin{verbatim}
>> prob = set(prob, 'Property1', Argument1, ..., 'Property2', ...);
\end{verbatim}

%In the next table the problem data are listed on the left column, while the corresponding property names are listed
%on the right column.
%\begin{center}
% use packages: array
%\begin{tabular}{|l|l|}
%\hline $x$ & \texttt{state} \\ 
%\hline $u$ & \texttt{input} \\ 
%\hline $t$ & \texttt{time} \\ 
%\hline $T$ & \texttt{horizon} \\ 
%\hline $f(t,x,u)$ & \texttt{dynamics} \\ 
%\hline $h(t,x,u)$ & \texttt{intcost} \\ 
%\hline $H(x)$ & \texttt{fcost} \\ 
%\hline $\CC_I$ & \texttt{iconstraint} \\ 
%\hline $\CC_T$ & \texttt{tconstraint} \\ 
%\hline $\CC_F$ & \texttt{fconstraint}, \texttt{fdirac} \\ 
%\hline $\mu_I(x_0)$ & \texttt{idirac}, \texttt{iuniform} \\
%\hline & \texttt{testtime}, \texttt{intconstraint} \\ \hline
%\end{tabular}
%\end{center}

To help the user remember the name of the properties the following convention has been used:
``\texttt{i}'' refers to the initial condition, ``\texttt{f}'' to the final condition,
``\texttt{t}'' to the trajectory, and ``\texttt{s}'' (for sum) to an integral.

In the sequel we list and describe all the properties that can be specified and their arguments.

\subsubsection{Setting the state variables}
\textbf{Property name:} \texttt{state}.
\\
\textbf{Number of arguments:} 1.
\\
\textbf{Default value:} empty.
\\
\textbf{Description:} Specifies state variables, a vector of
GloptiPoly class \texttt{mpol}.
\\
\textbf{Example}
\\
In this example we declare a variable \texttt{x} of dimension 2 and, after creating a POCP object \texttt{prob}, we set \texttt{x} as the state variable.
\begin{verbatim}
>> mpol x 2;
>> prob = pocp;
>> prob = set(prob, 'state', x);
\end{verbatim}
State variables can also have different names as in the following example
\begin{verbatim}
>> mpol x1 x2;
>> prob = pocp;
>> prob = set(prob, 'state', [x1; x2]);
\end{verbatim}
For more information on declaring GloptiPoly variables type
\begin{verbatim}
>> help mpol
\end{verbatim}
or check \cite{HenLasLof2007}.

\subsubsection{Setting the input variables}
\textbf{Property name:} \texttt{input}.
\\
\textbf{Number of arguments:} 1.
\\
\textbf{Default:} empty.
\\
\textbf{Description:} Specifies input variables, a vector of
GloptiPoly class \texttt{mpol}.
\\
\textbf{Example}
\\
Suppose a POCP object \texttt{prob} has been already defined. With the following commands a scalar \texttt{mpol} variable \texttt{u} is defined and then set as the input variable.
\begin{verbatim}
>> mpol u;
>> prob = set(prob, 'input', u);
\end{verbatim}

\subsubsection{Setting the time variable}
\textbf{Property name:} \texttt{time}.
\\
\textbf{Number of arguments:} 1.
\\
\textbf{Default:} empty.
\\
\textbf{Description:} Specifies the time variable,
a scalar of GloptiPoly class \texttt{mpol}. 
When the time variable is not specified but is needed internally by POCP, then it is defined automatically.
\\
\textbf{Example}
\\
Suppose a POCP object \texttt{prob} has already been defined. With the following commands a scalar \texttt{mpol}
variable \texttt{t} is defined and then set as the time variable.
\begin{verbatim}
>> mpol t;
>> prob = set(prob, 'time', t);
\end{verbatim}

\subsubsection{Setting the dynamics}
\textbf{Property name:} \texttt{dynamics}.
\\
\textbf{Number of arguments:} 1.
\\
\textbf{Default:} empty.
\\
\textbf{Description:} Specifies the system dynamics, a vector of
GloptiPoly class \texttt{mpol} of the same size as the state vector.
\\
\textbf{Example}
\\
Assume a POCP object \texttt{prob} has been created and the state vector \texttt{x} and input variable \texttt{u} have been set. In this example we set the linear dynamics specified by the state matrix \texttt{A} and input matrix \texttt{B}.
\begin{verbatim}
>> A = [0 1; 1 1];
>> B = [0; 1];
>> prob = set(prob, 'dynamics', A*x+B*u);
\end{verbatim}

\subsubsection{Setting the horizon length}
\textbf{Property name:} \texttt{horizon}.
\\
\textbf{Number of arguments:} 1.
\\
\textbf{Default:} 0.
\\
\textbf{Description:} Specifies the horizon length, a non-negative number. If set to 0, the horizon is free.
\\
\textbf{Example}
\\
Assume the POCP object \texttt{prob} has already been created. With the following command the horizon length is set to 1.
\begin{verbatim}
>> prob = set(prob, 'horizon', 1);
\end{verbatim}

\subsubsection{Setting the final cost} 
\textbf{Property name:} \texttt{fcost}.
\\
\textbf{Number of arguments:} 1.
\\
\textbf{Default:} 0.
\\
\textbf{Description:} Specifies the final cost, a polynomial of
GloptiPoly class \texttt{mpol} whose variables are elements of the state vector. If both the final cost and the integral cost are 0, then POCP minimizes the trace of the moment matrix of the occupation measure of the trajectory.
\\
\textbf{Example}
\\
Assume a POCP object \texttt{prob} has already been created and \texttt{x(1)} is one of the state variables. The following command sets the final cost to \texttt{x(1)\^{}2}.
\begin{verbatim}
>> prob = set(prob, 'fcost', x(1)^2);
\end{verbatim}

\subsubsection{Setting the integral cost} 
\textbf{Property name:} \texttt{scost}.
\\
\textbf{Number of arguments:} 1.
\\
\textbf{Default:} 0.
\\
\textbf{Description:} Specifies the integral cost, a polynomial of
GloptiPoly class \texttt{mpol} whose variables are the time variable and elements of the state and input vectors. For minimum time optimal control problems, set \texttt{scost} to \texttt{1} and set \texttt{horizon} to 0. If both the final cost and the integral cost are 0, then POCP minimizes the trace of the moment matrix of the occupation measure of the trajectory.
\\
\textbf{Example}
\\
Assume a POCP object \texttt{prob} has already been created and the state and input vectors are \texttt{x} and \texttt{u}, respectively. With the following command the integral cost is set to \texttt{x'*x+u'*u}.
\begin{verbatim}
>> prob = set(prob, 'scost', x'*x+u'*u);
\end{verbatim}
\subsubsection{Setting the dependence on time of the test functions}
\textbf{Property name:} \texttt{testtime}.
\\
\textbf{Number of arguments:} 1.
\\
\textbf{Default:} empty.
\\
\textbf{Description:} When the system dynamics does not depend on time and the horizon is free the test function does not depend on time. In the other cases the test function depends on time. However, this default behavior can be modified by setting the property \texttt{testtime}. When the argument is \texttt{true} (\texttt{false}) the test functions will (not) depend on time. To restore the default behavior, specify an empty (\texttt{[]}) argument.
\\
\textbf{Example}
\\
Assume a POCP object \texttt{prob} has already been defined. To make the test function depend also on time,
use the following command
\\
\begin{verbatim}
>> prob = set(prob, 'testtime', true);
\end{verbatim}

\subsubsection{Setting equality and inequality constraints on the initial condition}
\textbf{Property name:} \texttt{iconstraint}.
\\
\textbf{Number of arguments:} 1.
\\
\textbf{Default:} empty.
\\
\textbf{Description:} Specifices equality and inequality constraints on the initial condition,
a column vector of GloptiPoly class \texttt{supcon}.
\\
\textbf{Example}
\\
Assume the POCP object \texttt{prob} has already been created and \texttt{x(1)} and \texttt{x(2)} are the state variables. With the following command, both variables are constrained to be less than or equal to 1.
\begin{verbatim}
>> prob = set(prob, iconstraint, [x(1)<=1; x(2)<=1]);
\end{verbatim}
If the sum of \texttt{x(1)} and \texttt{x(2)} must be equal to 1, the following command should be used
\begin{verbatim}
>> prob = set(prob, iconstraint, x(1)+x(2)==1);
\end{verbatim}

\subsubsection{Setting the initial condition to a Dirac delta}
\textbf{Property name:} \texttt{idirac}.
\\
\textbf{Number of arguments:} 2 or 3.
\\
\textbf{Default:} unspecified, considered as a problem unknown.
\\
\textbf{Description:} This property should be set when one or more state variables take their value at time 0 accordingly to a distribution given by the sum of Dirac deltas. Two or three arguments can be specified. The first argument is a column vector containing the state variables involved. The second is a matrix whose columns contain the value of the variables (one column for each Dirac delta). The third argument (optional when the second argument has only one column) is a row vector specifying the probability of each Dirac delta.
\\
\textbf{Example}
\\
Assume a POCP object \texttt{prob} has already been created and \texttt{x(1)} and \texttt{x(2)} are the state variables. With the following command we set the initial condition to $x(0)=(1,1)$
\begin{verbatim}
>> prob = set(prob, 'idirac', [x(1); x(2)], [1; 1]);
\end{verbatim}
If the value at time 0 is $x(0)=(0, 1)$ with probability $0.8$ and $x(0)=(1, 1)$ with probability $0.2$, we should use
\begin{verbatim}
>> prob = set(prob, 'idirac', x, [0 1; 1 1]; [0.8 0.2]);
\end{verbatim}

\subsubsection{Setting the initial condition to a uniform distribution on a box}
\textbf{Property name:} \texttt{iuniform}.
\\
\textbf{Number of arguments:} 2.
\\
\textbf{Default:} unspecified, considered as a problem unknown.
\\
\textbf{Description:} This property should be set when one or more state varibles take their value at time 0 accordingly to a uniform distribution on a box. Two arguments must be specified. The first argument is a column vector containing the state variables involved. The second is a matrix with two columns and the number of rows equal to the length of the first argument, specifying the intervals.
\\
\textbf{Example}
\\
Assume a POCP object \texttt{prob} has already been created and \texttt{x(1)} and \texttt{x(2)} are the state variables.
With the following command we set the value of \texttt{x(1)} to be uniformly distributed on the interval $[-1,~1]$ and the value of \texttt{x(2)} on the interval $[0,~2]$
\begin{verbatim}
>> prob = set(prob, 'iuniform', [x(1); x(2)], [-1 1; 0 2]);
\end{verbatim}
The properties \texttt{idirac} and \texttt{iuniform} can also be used in the same problem, like in the following
\begin{verbatim}
>> prob = set(prob, 'iuniform', x(1), [-1 1]);
>> prob = set(prob, 'idirac', x(2), -1);
\end{verbatim}

\subsubsection{Setting equality and inequality constraints on the final condition}
\textbf{Property name:} \texttt{fconstraint}.
\\
\textbf{Number of arguments:} 1.
\\
\textbf{Default:} empty.
\\
\textbf{Description:} Specifies equality and inequality constraints on the state variables at the end of the horizon. 
The usage is the same as for the property \texttt{iconstraint}.

\subsubsection{Setting the final condition to a Dirac delta}
\textbf{Property name:} \texttt{fdirac}.
\\
\textbf{Number of arguments:} 2 or 3.
\\
\textbf{Default:} unspecified, considered as a problem unknown.
\\
\textbf{Description:} This property should be set when one or more state variables take their value at the end of the horizon accordingly to a distribution given by the sum of Dirac deltas.
The usage is the same as for the property \texttt{idirac}.

\subsubsection{Setting equality and inequality constraints on the variables along the trajectory}
\textbf{Property name:} \texttt{tconstraint}
\\
\textbf{Number of arguments:} 1.
\\
\textbf{Default:} empty.
\\
\textbf{Description:} Specifies equality and inequality constraints on the system variables
along the trajectory, a column vector of GloptiPoly class \texttt{supcon}.
\\
\textbf{Example}
\\
Assume a POCP object \texttt{prob} has already been created and \texttt{x} is the scalar state, \texttt{u} is the input, and \texttt{t} is the time variable. With the following command we set the constraints \texttt{x>=0} and \texttt{x+u+t<=1}
\begin{verbatim}
>> prob = set(prob, 'tconstraint', [x>=0; x+u+t<=1]);
\end{verbatim}

\subsubsection{Setting integral constraints}
\textbf{Property name:} \texttt{sconstraint}.
\\
\textbf{Number of arguments:} 1.
\\
\textbf{Default:} empty.
\\
\textbf{Description:} Specifices constraints on the integral of some variables, a column vector of
GloptiPoly class \texttt{momcon}.
\\
\textbf{Example}
\\
Assume a POCP object \texttt{prob} has already been created and \texttt{u} is the input variable.
To enforce the integral constraint $\int_0^T u(t)^2 dt \leq 1$, use the following command
\begin{verbatim}
>> prob = set(prob, 'sconstraint', mom(u^2)<=1);
\end{verbatim}

\subsubsection{Some remarks on the usage of the \texttt{set} command}
\begin{itemize}
\item The system variables (state vector, input vector, and time) must be  set before they are used. Constraints on the initial state cannot be entered if the state vector has not been set.
When more properties are specified with the same \texttt{set} command, the system variables should be specified first.
\item Once the system variables (state vector, input vector, and time) have been specified, they cannot be modified any further.
\item The current version of the software does not check for conflicts between different kinds of constraints. E.g., if the value at time $0$ for a state variable is assigned by setting the property \texttt{idirac} and then the same variable appears in the constraints set by \texttt{iconstraint}, the behavior of POCP is unpredictable.
\end{itemize}

\subsection{The \texttt{solvepocp} command}\label{su:solvepocp}
\textbf{Short description:} Builds and solves an SDP relaxation of a POCP.
\\
\textbf{Syntax and usage:} When used in the form
\begin{verbatim}
>> [status, cost, measures] = solvepocp(prob, degree);
\end{verbatim}
or, equivalently,
\begin{verbatim}
>> [status, cost, measures] = solvepocp(prob, 'mom', degree);
\end{verbatim}
the command
\texttt{solvepocp} builds and solves an SDP relaxation where moments of degree up to \texttt{degree} are used. The input
argument \texttt{'mom'} stands for moments. The output parameters have the following meaning:
\begin{itemize}
\item When \texttt{status} is strictly negative, the SDP problem
is infeasible or could not be solved. Otherwise,
the SDP problem could be solved;
\item \texttt{cost} is a lower bound on the objective function of the POCP
if \texttt{status} is non-negative;
\item \texttt{measures} is a structure containing the measures used by GloptiPoly to solve the relaxation of the POCP.
Once the relaxation is solved, \texttt{measures} contains fields \texttt{measures.initial}, \texttt{measures.final}, and \texttt{measures.trajectory}. If the initial or final
conditions have been assigned, the corresponding measures are empty (\texttt{[]}).
To extract the moment matrix corresponding to one of the measures, use the GloptiPoly function \texttt{mmat} 
in combination with \texttt{double}. For example, to retrieve the moment matrix corresponding to the trajectory occupation measure, use the following command
\begin{verbatim}
>> double(mmat(measures.trajectory));
\end{verbatim}
\end{itemize}
It is also possible to call \texttt{solvepocp} by specifying the degree of the test function to be used when defining the moment problem:
\begin{verbatim}
>> [status, cost, measures] = solvepocp(prob, 'tf', degree);
\end{verbatim}
The input argument \texttt{'tf'} stands for test function.
As a rule of thumb, \texttt{degree} should be at least twice the degree of the dynamics.
Increasing the value of \texttt{degree} improves the quality of the approximation,
but it also increases the size of the SDP problem solved.

The command \texttt{solvepocp} can also be called by specifying four output parameters:
\begin{verbatim}
>> [status, cost, measures, vf] = solvepocp(prob, degree);
\end{verbatim}
Using this syntax, a subsolution of the HJB equation is calculated and stored in \texttt{vf} (for value function). This approximation is a polynomial of GloptiPoly class \texttt{mpol}.
If the option \texttt{'tf'} is used, \texttt{degree} corresponds to the degree of the subsolution of the
HJB equation. The subsolution \texttt{vf} is computed by using the
YALMIP toolbox, which must be installed properly.

\section{First example}
In this section we show how POCP can be used to design a suboptimal controller.
Consider as in \cite{HenLasSav2008}
the nonlinear system
\[
\begin{bmatrix} \dot x_1(t) \\ \dot x_2(t) \end{bmatrix} =
\begin{bmatrix} x_2(t) - x_1(t)^3 + x_1(t)^2 \\ u(t) \end{bmatrix}.
\]
The optimal control problem consists in driving to the origin the initial states from
the set $(x_1,x_2)\in \SS =[-1,1]\times[-1,1]$ by minimizing the cost functional $\int_0^Th(x(t),u(t))dt$ with
\[
h(x(t), u(t))=x_1(t)^2+x_2(t)^2+\frac{u(t)^2}{100}.
\]
The terminal time $T$ is free.

To solve this problem and design a control law, we can define a POCP problem and use the command \texttt{solvepocp} 
with four output arguments. By setting the initial condition to be uniformly distributed on the box $\SS$, the fourth output argument of \texttt{solvepocp}
contains a polynomial subsolution of the HJB equation which approximates the value function along all the optimal trajectories starting from $\SS$, see \cite{HenLasSav2008}.
Therefore, a natural choice to derive a control law consists in solving the following minimization problem
\begin{equation}\label{eq:hjb}
\min_{u} \left[\frac{\partial V_f(x,t)}{\partial t} + \nabla_x V_f(x)f(x,u) + h(x,u)\right],
\end{equation}
where $V_f$ represents the approximated value function.
Since the dynamics is affine in $u$ and $h(\dot)$ is quadratic in $u$, the global minimizer of (\ref{eq:hjb}) can be calculated by using first order optimality conditions:
\begin{equation}\nonumber
u(x):=-50\nabla_x \varphi(x) \begin{bmatrix}0 & 1\end{bmatrix}^T.
\end{equation}
In the sequel we show a Matlab POCP script which calculates a polynomial approximation of the value function \texttt{vf} and defines the control law \texttt{u\_x}.
\begin{verbatim}
% variable declaration
mpol x 2
mpol u

% problem definition
prob = pocp( ...
       'state', x, ...
       'input', u, ...
       'dynamics', [x(2)+x(1)^2-x(1)^3; u], ...
       'horizon', 0, ...
       'iuniform', x, [-1 1; -1 1], ...
       'fd', x, [0;0], ...
       'scost', x'*x+u^2/100);

% problem solved with test function up to degree 8
[status, cost, measures, vf] = solvepocp(prob, 'tf', 8);

% control law definition
u_x = -50 * diff(vf, x)*[0; 1];
\end{verbatim}
The value of \texttt{u\_x} can be easily evaluated by using the
GloptiPoly commands \texttt{assign} and \texttt{double}.
For example, to obtain a value for the control at $x=(0.5, ~0.5)$, the following commands can be used.
\begin{verbatim}
>> assign([x; u], [0.5 0.5 0]');
>> controlvalue = double(u_x);
\end{verbatim}
Notice that when we assign the value of \texttt{x} we must assign also the value of \texttt{u} (any value can be used). For more information on the commands \texttt{assign} and \texttt{double}, see the documentation of GloptiPoly.

\begin{figure}[htb]
\centering{
\includegraphics[width=0.7\textwidth]{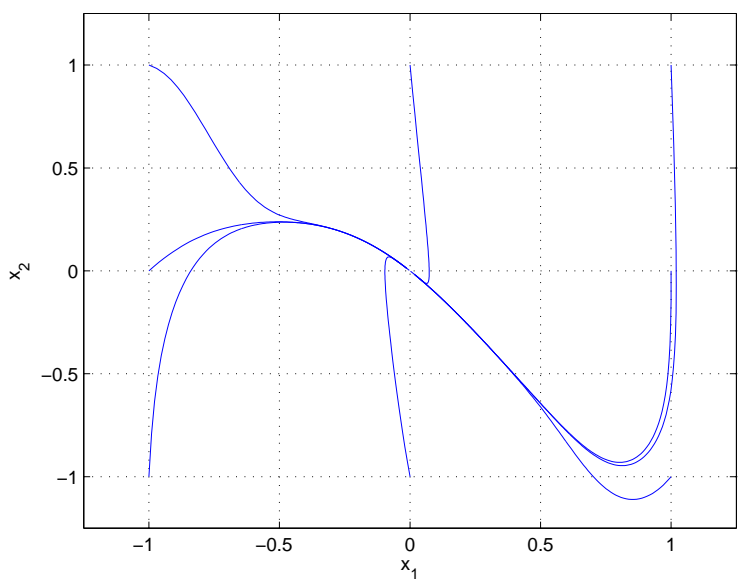} }
\caption{Trajectories obtained by applying the control law \texttt{u\_x}.
\label{fig:nonlinear}}
\end{figure}

Figure \ref{fig:nonlinear} shows some trajectories starting from $\SS$ and obtained by simulating the system with the calculated control law.

\section{Second example}
This example is taken from \cite[Example 3.27]{SepJanKok1997}.
For the nonlinear system
\[
\begin{bmatrix} \dot x_1(t) \\ \dot x_2(t) \end{bmatrix} =
\begin{bmatrix} -x_1(t)^3 + x_1(t)u(t) \\ u(t) \end{bmatrix}
\]
and the cost functional $\int_0^{\infty} (x_2(t)^2+u(t)^2)dt$,
the solution of the HJB equation is $x_2(t)^2$, which is only positive
semidefinite. Since the solution is polynomial, we can show with this
example that the solution can be found exactly (and only using test 
functions up to degree 2):
\begin{verbatim}
mpol x 2
mpol u

p=pocp();
p=set(p,'state',x);
p=set(p,'input',u);
p=set(p,'dynamics',[-x(1)^3+x(1)*u; u]);
p=set(p,'scost',x(2)^2+u^2);
p=set(p,'idirac',x, [1; 1]);
p=set(p,'fdirac',x, [0; 0]);
p=set(p,'tconst', [x>=-1.1; x<=1.1]);

[status,cost,measures,vf] = ...
 solvepocp(p,'tf', 2);

disp('Subsolution of HJB equation=');
vf
\end{verbatim}
Running the above script we obtain the following output:
\begin{verbatim}
Scalar polynomial
x(2)^2
\end{verbatim}

\section{Variable scaling}
To improve the numerical behavior of the SDP solver used by POCP, the variables (state, input and time) could be scaled. Scaling can be achieved by using the GloptiPoly command \texttt{scale}. For example, to substitute the variable \texttt{x} with \texttt{2*x}, use the following commands
\begin{verbatim}
>> mpol x
>> scale(x, 2)
\end{verbatim}
Scaling changes only the internal representation of the variables. From the user point of view nothing changes,
and there is no need to scale the result of the optimization problem.
For more information about the \texttt{scale} command, see the documentation of GloptiPoly.
If possible, all the variables should be scaled within the interval $[-1, +1]$.

\section{Some ideas for new features}
We end this guide by listing ideas for improvements and new features of POCP:
\begin{itemize}
% \item It is possible to extend the occupation measure approach to discrete-time and hybrid systems. Following this idea, POCP could also be extended to deal with these systems.

\item To further extend the class of continuous-time problems considered, more flexibility on the initial and final time (currently $0$ and $T$) could be permitted. If the cost
\begin{equation}
J(t_i,t_f,x(0),u(\cdot))=\int_{t_i}^{t_f} h(t, x(t), u(t)) dt + H(x(t_f)).
\end{equation}
is considered, the problem
\begin{equation}
\begin{array}{ll}
\displaystyle\min_{u(\cdot)} & \int J(t_i,T,x_i,u(x_i,t)) d\mu_i(t_i,x_i) \\
& (t,x(t),u(t))\in \CC_{T} \\
& x(T)\in \CC_{F}
\end{array}
\end{equation}
could be modeled with POCP, with $\mu_i$ an assigned probability measure on $\R\times\R^n$.
Setting $\mu_i$ to be a uniform measure on the set $\{(t,x):t\in[0,1],~ x_j\in[-1,-1]~\forall j=1,\dots,n\}$, we could obtain a subsolution of the HJB equation
which gives a good approximation of the value function on a box centered at the origin of the state space and for all the values of $t$ in the interval $[0,1]$, see \cite{HenLasSav2008}.

\item Version 3.3 of GloptiPoly cannot be used to retrieve information on the dual variables of the moment problem. As a workaround to calculate a subsolution of the HJB equation, POCP currently uses GloptiPoly in combination with YALMIP. If GloptiPoly will include this feature in some future version, POCP will not require YALMIP.

\item GloptiPoly does not support currently the assignment of a measure (e.g., it cannot be specified that a variable is uniformly distributed on an interval). If GloptiPoly will include this feature in some future version, the POCP code
could be significantly simplified.
\end{itemize}

\end{document}